\documentclass[11pt]{article}
\usepackage{latexsym}
\usepackage{maa-monthly}

\usepackage{amsfonts, amssymb, amsthm, amsmath, hyperref, graphicx}
\usepackage[mathscr]{eucal}

\newtheorem{theorem}{Theorem}
\newtheorem*{CarlitzThm}{Carlitz's Theorem for Class Number 2}
\newtheorem*{CarlitzThmR}{Carlitz's Theorem Redux}
\newtheorem*{Nutshell}{Class Number 2 in a Nutshell}
\newtheorem*{Fundamental}{The Fundamental Theorem of Ideal Theory}
\newtheorem{corollary}[theorem]{Corollary}

\newtheorem{proposition}[theorem]{Proposition}

\newtheorem{lemma}[theorem]{Lemma}

\theoremstyle{definition}
\newtheorem{example}[theorem]{Example}

\newtheorem{definition}[theorem]{Definition}

\newcommand{\mf}{\mathfrak}

\newcommand{\dis}{\mathbf{d}}

\begin{document}
\title{So What is Class Number 2?}
 
%%
%% Now edit the following to give your name and address:
%% 

\markright{Class Number 2}
%\author{Anonymous}
\author{Scott T. Chapman}

\maketitle

\begin{abstract}  Using factorization properties, we give several characterizations for a ring of algebraic integers to have class number at most 2.
\end{abstract}
\medskip

\section{Introduction.}

I was recently at an algebra colloquium when some questions involving small class numbers of algebraic number rings arose.  Of the 30 or so participants, almost everyone in the room recognized that an algebraic number
ring is a unique factorization domain (or UFD) if and only if its class number is one (i.e., the ideal class group of $R$ is trivial).  Almost no one in the room was aware of the following theorem of Carlitz, which is well known among mathematicians
who work in the theory of nonunique factorizations (see \cite{GHKb} for a general reference to this area).  

\begin{CarlitzThm}\cite{Ca}\label{Carlitz}
Let $R$ be an algebraic number ring.  $R$ has class number at most 2 if and only if whenever $\alpha_1,\ldots, \alpha_n,\beta_1, \ldots, \beta_m$ are irreducible elements of $R$
with 
\begin{equation}\label{CarlitzEq}\tag{$\dagger$}
\alpha_1\cdots \alpha_n=\beta_1\cdots \beta_m
\end{equation}
then $n=m$.
\end{CarlitzThm}

An integral domain $D$ in which each nonzero nonunit can be factored as a product of irreducible elements is known as an \textit{atomic domain}.   An atomic domain $D$ that satisfies the condition in Carlitz's theorem (i.e., satisfies \eqref{CarlitzEq}) is
known as a \textit{half-factorial domain} (or \textit{HFD}).  Notice that a UFD is an HFD and hence, if $R$ exactly has class number 2, it is an example of an HFD that is not a UFD (the classic such example is $\mathbb{Z}[\sqrt{-5}]$ and the nonunique factorization $6=2\cdot 3 = (1+\sqrt{-5})(1-\sqrt{-5})$).  Thus, the Carlitz theorem can be restated as follows.

\begin{CarlitzThmR}\label{CarlitzR}
  Let $R$ be an algebraic number ring.  $R$ has class number at most 2 if and only if $R$ is a half-factorial domain.
\end{CarlitzThmR}

Carlitz's theorem was the beginning of quantitative and qualitative research into nonunique factorizations in integral domains and monoids.  This research began with papers concerning HFD's (see \cite{Sk}, \cite{Z1}, \cite{Z2}, and a comprehensive survey article  \cite{ChapCoy}) and has expanded into the study of a host of combinatorial constants that measure deviation of factorizations from the UFD condition.  The purpose of this article is not to deeply explore the general topic of factorization, but to give a series of factorization-inspired characterizations of class number 2.  We do this solely in terms of algebraic number fields and thus avoid the abstraction and generality that more difficult factorization problems entail.  Our characterizations will involve constants of increasing complexity, 
and in light of this, we will offer the various needed definitions directly before each result.  We hope that our work gives the reader a better appreciation of Carlitz's theorem and its related substantive factorization problems.  For those who want a more in-depth treatment of nonunique factorizations, several recent papers on this topic can be found in this \textsc{Monthly} (\cite{ChapBag}, \cite{Ger1}, \cite{OnPel}).

%%%%%%%%%%%%%
%%% Revised introduction
%%%%%%%%%%%%%%%%%%%%%%%%%%

Throughout we assume an understanding of abstract algebra at the level of \cite{Ga} and a basic familiarity with algebraic number theory at the level of \cite{Ma}.  (An approach that might be more friendly to a novice can be found in \cite{PD}.)  
For clarity, we review the basic definitions necessary for the remainder of this work.  If $\mathbb{Q}$ represents the field of rational numbers, then an algebraic number field $K$ is any finite extension of $\mathbb{Q}$.  An element $\alpha \in K$ is an algebraic integer if it is a root of a monic polynomial in $\mathbb{Z}[X]$.  By \cite[Theorem 6.2]{PD}, the set $R$ of algebraic integers in $K$ is an integral domain, which we refer to as an algebraic number ring.   

When dealing with an algebraic number ring $R$, we use the usual notions of divisibility from the theory of integral domains.  Let $\mathcal{A}(R)$ represent the set of irreducible elements (or atoms) of $R$, $\mathcal{U}(R)$ the set of units of $R$, and $R^\bullet$ the set of nonzero nonunits of $R$.  Recall that $x$ and $y$ in $R$ are associates if there is a unit $u\in R$ with $x=uy$.  If $x, y$, and $z$ are in $R$ with $y=xz$, then we say that $x$ divides $y$ and denote this by $x\, |\, y$. 

Let $\mathcal{I}(R)$ denote the set of ideals of $R$.  If $x\in R$, then let $(x)$ represent the principal ideal generated by $x$ and $\mathcal{P}(R)$ the subset of $\mathcal{I}(R)$ consisting of principal ideals of $R$.  For $I$ and $J$ in $\mathcal{I}(R)$, set
\[
IJ = \left\{\sum_{i=1}^n a_ib_j \,\mid\, a_i\in I \mbox{  and  }b_j\in J\right\}.
\]
Using \cite[Theorem 8.1]{PD}, it is easy to argue that $IJ$ is another ideal of $R$ which is known as the product of 
$I$ and $J$.  If $I, J$, and $K$ are ideals
of $R$ with $J=IK$, then we borrow the notation used above for elements and say that $I \, |\, J$.  

Define an equivalence relation on $\mathcal{I}(R)$ by $I\sim J$ if and only if there exist $\alpha$ and $\beta$ in $R$ with $(\alpha)I=(\beta)J$.  If $[I]$ represents the equivalence class of the ideal $I$ under $\sim$, then by \cite[Lemma 10.1]{PD} the operation
\[
[I]+[J] = [IJ]
\]
is well-defined.  By \cite[Theorem 8.13]{PD}, the set $\mathcal{C}(R)=\mathcal{I}(R)/\sim$ forms an abelian group under $+$ called the class group of $R$.  By \cite[Theorem 10.3]{PD}, $\lvert \mathcal{C}(R) \rvert$
is finite and is known as the class number of $R$.  As previously mentioned, classical algebraic number theory (\cite[Theorem 9.4]{PD}) asserts that $R$ is a unique factorization domain if and only if its class number is one.  
Throughtout the rest of our work we will use freely the fact asserted in \cite{Ca} that every ideal class of $\mathcal{C}(R)$ contains infinitely many prime ideals.

To completely understand how elements factor in an algebraic number ring $R$, we will need this fundamental result concerning the factorizations of ideals in $R$.
\begin{Fundamental} \cite[Theorem 8.27]{PD}
Let $R$ be an algebraic number ring.  If $I$ is an ideal of $R$, then there exists a unique (up to order) list of not necessarily distinct prime ideals
$\mf{p}_1$, $\mf{p}_2, \ldots , \mf{p}_k$ of $R$ such that 
\begin{equation}\label{fund}\tag{$\star$}
I=\mf{p}_1\mf{p}_2\cdots \mf{p}_k.
\end{equation}
\end{Fundamental}
\noindent The key to comprehending factorizations in $R$ lies in understanding products of the form \eqref{fund} where $\sum_{i=1}^k [\mf{p}_i]=0$ in $\mathcal{C}(R)$ (see Lemma \ref{irreducibles} below).

%%%%%%%%%%%%%%%%%%%%%%%%%%%%%%%%%%
%%%%%%%%%%%%%%%%%%%%%%%%%%%%%%%%%%
%%%%%%%%%%%%%%%%%%%%%%%%%%%%%%%%%%

\section{More on the Carlitz Characterization.}  

We open with a few simple lemmas which will prove useful, especially in our later work.  The first will characterize the irreducible elements of $R$ in terms of the class group.

\begin{lemma}\label{irreducibles}
Let $R$ be an algebraic number ring and $x$ a nonzero nonunit of $R$ with
\[
(x)=\mf{p}_1\cdots \mf{p_n},
\]
where $n\geq 1$ and the $\mf{p}_i$'s are not necessarily distinct prime ideals of $R$. The element $x$ is irreducible in $R$ if and only if 
\begin{enumerate}
\item $\sum [\mf{p}_i] =0$, and
\item if $S\subsetneq \{1,\ldots n\}$ is a nonempty subset then $\sum_{i \in S}[\mf{p}_i]\neq 0$.
\end{enumerate}
\end{lemma}

\begin{proof}  ($\Rightarrow$) That $\sum [\mf{p}_i] =0$ follows from the definition of the class group.  Suppose there is a proper subset $S$ of $ \{1,\ldots n\}$ with $\sum_{i \in S}[\mf{p}_i]= 0$.  Let $S'=\{1,\ldots ,n\}-S$.  Then both
\[
\sum_{i \in S}[\mf{p}_i]=0\mbox{  and  } \sum_{i \in S'}[\mf{p}_i]=0
\]
and hence there are nonunits $y$ and $z$ in $R$ with 
\[
(y)=\prod_{i\in S}\mf{p}_i \mbox{  and  } (z)= \prod_{i\in S'}\mf{p}_i.
\]
Thus, there is a unit $u\in R$ with $x=uyz$ and $x$ is not irreducible.

($\Leftarrow$)  Suppose that $x=yz$ in $R$.  By the fundamental theorem of ideal theory in $R$, there are nonempty subsets $S, S'$ of $\{1, \ldots ,n\}$ so that
\[
(y)=\prod_{i\in S}\mf{p}_i\mbox{  and  } (z)=\prod_{i\in S'} \mf{p}_i.
\]
Then $\sum_{i\in S} [\mf{p}_i] =0$ contradicting condition (2).  This completes the proof.
\end{proof}

\begin{example}\label{ex1}  We illustrate the results of the lemma with some examples.  Let $\mf{p}$ be a nonprincipal prime ideal of $R$ with $|[\mf{p}]|=n$ (where $|[\mf{p}]|$ represents of order of $[\mf{p}]$ in $\mathcal{C}(R)$).  Then
\[
\mf{p}^n=(x),
\]
where $x$ is irreducible in $R$.  Moreover, if $\mf{q}$ is any prime ideal taken from class $-[\mf{p}]$, then 
\[
\mf{p}\mf{q}=(y),
\]
where $y$ is irreducible in $R$.  Hence, in the case where $|\mathcal{C}(R)|=2$, an irreducible element takes one of three forms: 
\begin{enumerate}
\item[(i)] $\alpha$ where $(\alpha)=\mf{p}$ for a principal prime ideal $\mf{p}$ of $R$;
\item[(ii)] $\alpha$ where $(\alpha)=\mf{p}^2$ for a nonprincipal prime ideal $\mf{p}$ of $R$;
\item[(iii)] $\alpha$ where $(\alpha)=\mf{p}\mf{q}$  where $\mf{p}$ and $\mf{q}$ are distinct nonprincipal prime ideals of $R$.
\end{enumerate}
In case (i), the irreducible $\alpha$ is actually a prime element; in case (ii), $\alpha$ is called \textit{ramified}; and in case (iii), $\alpha$ is called \textit{split}.
\end{example}

Lemma \ref{irreducibles} implies some important finiteness conditions.  A sequence of elements $g_1, \ldots ,g_n$ from an abelian group $G$ that satisfies the sum condition in the lemma (i.e., $g_1+\cdots + g_n=0$ and no proper subsum of this sum is zero)
is known as a \textit{minimal zero-sequence}.  A good reference on the interplay between factorizations in an algebraic number ring and minimal zero-sequences is \cite{GRu}. An elementary exercise (see for example \cite{Chap}) shows that the number of minimal zero-sequences in a finite abelian group is finite.  Since there are finitely many, there is a finite constant known as $D(G)$ that bounds above the number of elements in this minimal zero-sequence.  The computation of $D(G)$, known as the Davenport constant of $G$, is elusive and better left to our references (\cite{Chap} is a good source).   These facts imply the following corollary.

\begin{corollary}\label{thecor}
Let $x\in R^\bullet$ where $R$ is an algebraic number ring.  
\begin{enumerate}
\item[(1)] The element $x$ has finitely many nonasscociated irreducible factorizations.
\item[(2)] If $x$ is irreducible and $(x) = \mf{p}_1\cdots \mf{p}_k$, then $k\leq D(\mathcal{C}(R))$.  
\end{enumerate}
\end{corollary}

Having established that irreducible factorizations are essentially finite in number, we produce one below which will be of particular interest.

\begin{lemma}\label{carlitzlemma}
Let $R$ be a ring of algebraic integers of class number greater than 2.  Then there are not necessarily distinct irreducible elements $\alpha_1, \alpha_2, \beta_1, \beta_2$, and $\beta_3$ such that
\begin{equation}\label{basic}\tag{$\ddagger$}
\alpha_1\alpha_2 = \beta_1\beta_2\beta_3.
\end{equation}
\end{lemma}

\begin{proof}
Suppose that $\mathcal{C}(R)$ contains an element $g$ with $|g|=n>2$.  Let $\mf{p}_1$ be a prime ideal of $R$ taken from class $g$, $\mf{p}_2$ a prime ideal taken from class $2g$, $\mf{p}_3$ a prime ideal taken 
from class $(n-2)g$, and $\mf{p}_4$ a prime ideal taken from class $(n-1)g$.  (In the cases $n=3$ or $4$, you can pick these ideals distinctly.)  Define the irreducible elements $\alpha$, $\beta$, $\gamma$, and $\delta$ of $R$ by
\begin{enumerate}
\item $(\alpha) = \mf{p}_1\mf{p}_4$,
\item $(\beta) = \mf{p}_1^2\mf{p}_3$,
\item $(\gamma) = \mf{p}_2\mf{p}_3$,
\item $(\delta) = \mf{p}_2\mf{p}_4^2$.
\end{enumerate}
The ideal equation $(\mf{p}_1^2\mf{p}_3)(\mf{p}_4^2\mf{p}_2) = (\mf{p}_1\mf{p}_4)^2(\mf{p}_2\mf{p}_3)$ yields that 
\[
\beta\delta = u\alpha^2\gamma
\]
for some $u\in \mathcal{U}(R)$.

If all the nonidentity elements of $\mathcal{C}(R)$ are of order 2, then let $g_1$ and $g_2$ be such elements with $g_1\neq g_2$.  Suppose further that $g_3=g_1+g_2$.  Thus, $g_1, g_2$, and $g_3$ are all
distinct elements of $\mathcal{C}(R)$ of order 2.  If $\mf{p}_1$, $\mf{p}_2$, and $\mf{p}_3$ are prime ideals of $\mathcal{C}(R)$ taken from the classes $g_1$, $g_2$, and $g_3$ respectively, then
\[
\mf{p}_1^2 = (\beta_1),\; \mf{p}_2^2 = (\beta_2),\; \mf{p}_3^2=(\beta_3),\; \mbox{ and } \mf{p}_1\mf{p_2}\mf{p_3} = (\alpha)
\]
with $\beta_1$, $\beta_2$, $\beta_3$, and $\alpha$ irreducible elements of $R$.  Thus in $R$ we have
\[
\alpha^2 = u\beta_1\beta_2\beta_3
\]
for some unit $u$ of $R$.  This completes the proof.
\end{proof}

We are now in a position to offer a very short proof of Carlitz's theorem.

\begin{proof}[Proof of Carlitz's Theorem] ($\Leftarrow$)  If $R$ is half-factorial, then \eqref{basic} implies that $|\mathcal{C}(R)|\leq 2$. 

($\Rightarrow$)  Let $x\in R^\bullet$ with 
\[
(x)=\mf{q}_1\cdots \mf{q}_n\mf{p}_1\cdots \mf{p}_m,
\]
where the prime ideals $\mf{q}_i$ are principal and the prime ideals $\mf{p}_j$ are not principal.  By our remarks in Example \ref{ex1}, $m$ is even and any factorization of $x$ into irreducibles has length $n + \frac{m}{2}$.  Thus (2) holds and the proof 
is complete.
\end{proof}

\section{Characterizations Involving the Length Set.}  If $R$ is an algebraic number ring and $x$ a nonzero nonunit of $R$, then set
\[
\mathcal{L}(x)=\{ k\, |\, \exists \mbox{ irreducibles } \alpha_1, \ldots , \alpha_k \in R \mbox{ with  } x=\alpha_1\cdots \alpha_k\}.
\]
The set $\mathcal{L}(x)$ is known as the set of lengths of $x$ and a general \textsc{Monthly} survey on this topic can be found in \cite{Ger1}.  Corollary \ref{thecor} implies
that $|\mathcal{L}(x)|<\infty$ for any $x\in R^\bullet$.  By Carlitz's theorem, if $R$ has class number 2, then $\mathcal{L}(x)=\{k\}$ for some $k\in\mathbb{N}_0$, and if $|\mathcal{C}(R)|>2$, then Lemma \ref{carlitzlemma} implies that 
there is an $x\in R$ with $|\mathcal{L}(x)|>1$.    Set
\[
L(x)=\max\, \mathcal{L}(x), \; \ell(x)=\min\, \mathcal{L}(x),
\]
and 
\[
\rho(x)=\frac{L(x)}{\ell (x)}.
\]
Since $L(x)<\infty$, $\rho(x)$ is a rational $q\geq 1$ which is known as the \textit{elasticity} of $x$ in $R$.  We can turn this combinatorial constant into a global descriptor by setting
\[
\rho(R) = \sup\{\rho(x)\, |\, x\in R\}.
\]
Hence, $R$ is half-factorial if and only if $\rho(R)=1$ and by Lemma \ref{carlitzlemma}, if $R$ has class number greater than 2, then $\rho(R)\geq \frac{3}{2}$.   A 
detailed study of elasticity in number rings can be found in \cite{V} and a more general survey on the subject in \cite{And}.   In \cite{V} it is established that 
\[
\rho(R)= \frac{D(\mathcal{C}(R)}{2},
\]
where again $D(\mathcal{C}(R))$ represents Davenport's constant.

A more precise version of the elasticity has recently become popular in the literature.  Let $k\in \mathbb{N}$ and set
\[
\rho_k(R) =\sup\{\sup\, \mathcal{L}(x)\, |\, \min\, \mathcal{L}(x) \leq k\mbox{  for  } x\in R^\bullet\}.
\]
Using Corollary \ref{thecor} along with \cite[Proposition 1.4.2]{GHKb}, the fact that $R$ is an algebraic number ring yields a slightly simpler version of this definition:
\[
\rho_k(R) =\sup\{\max\, \mathcal{L}(x)\, |\, k \in \mathcal{L}(x) \mbox{  for  } x\in R^\bullet\}.
\]

We prove a few convenient facts concerning the $\rho_k(R)$'s.

\begin{lemma}\label{rhok}
If $R$ is an algebraic number ring, then the following assertions hold.
\begin{enumerate}
\item[(1)] $\rho_1(R)=1$.
\item[(2)] $\rho_k(R) \geq k$ for all $k\in \mathbb{N}$.
\item[(3)] For each $k\in \mathbb{N}$, $\rho_k(R)<\infty$.
\item[(4)] For each $k\in \mathbb{N}$, $\rho_k(R)<\rho_{k+1}(R)$.
\end{enumerate}
\end{lemma}

\begin{proof}
The proof of (1) follows directly from the definition of an irreducible element. For (2), if $x$ is a prime element of $R$, then $\mathcal{L}(x^k)=\{k\}$ so $k\in\mathcal{L}(x^k)$ and $k=\max \mathcal{L}(x^k)$.  That $\rho_k(R)\geq k$ now follows.

For (3), suppose that $k\in \mathcal{L}(x)$ for some $x\in R^\bullet$.  Thus $x=\alpha_1\cdots \alpha_k$ where each $\alpha_i\in \mathcal{A}(R)$.  Write each $(\alpha_i)=\mf{p}_{i,1}\cdots \mf{p}_{i,t_i}$ where each $\mf{p}_{i,j}$ is a prime
ideal of $R$.  By our previous comment, each $t_i\leq D(\mathcal{C}(R))$.  Thus, $(x)$ factors into at most $k\cdot D(\mathcal{C}(R))$ prime ideals, which also bounds the length of a factorization of $x$ into irreducibles.  Hence $\max \mathcal{L}(x)\leq
k\cdot D(\mathcal{C}(R))$ for each $x\in R^\bullet$ and thus $\rho_k(R)\leq k\cdot D(\mathcal{C}(R))$.   

For (4), suppose $m=\rho_k(R)$.  Then there are irreducible elements $\alpha_1,\ldots, \alpha_k$,  $\beta_1,\ldots, \beta_m$ of $R$ with $\alpha_1\cdots \alpha_k=\beta_1\cdots \beta_m$.  If $x$ is any irreducible element of $R$, then
$x\alpha_1\cdots \alpha_k=x\beta_1\cdots \beta_m$ and hence $\rho_{k+1}(R)\geq m+1 > \rho_k(R)$.
\end{proof}

The true relationship between $\rho(R)$ and the $\rho_k(R)$'s can be found in \cite[Proposition 6.3.1]{GHKb}:
\[
\rho(R) = \sup \left\{ \frac{\rho_k(R)}{k}\, |\, k\in \mathbb{N}\right\} = \lim_{k\rightarrow \infty} \frac{\rho_k(R)}{k}.
\]

Lemma \ref{carlitzlemma} again allows us to make an immediate deduction.  (Part of this result can be found prior to \cite[Proposition 1.4.2]{GHKb}.)

\begin{theorem}
Let $R$ be an algebraic number ring.  The following statements are equivalent.
\begin{enumerate}
\item[(1)] $R$ has class number less than or equal to 2.
\item[(2)] $\rho(R)=1$.
\item[(3)] $\rho_2(R) = 2$.
\item[(4)] $\rho_k(R) = k$ for some $k\geq 2$.
\item[(5)] For all irreducibles $x$ and $y$ in $R$, $\mathcal{L}(xy) = \{2\}$.
\item[(6)] For all irreducibles $x$ and $y$ in $R$, $|\mathcal{L}(xy)|=1$.
\end{enumerate}
\end{theorem}

\begin{proof}  Assertions (1) and (2) are equivalent by the Carlitz theorem.  If (2) holds, then every $\mathcal{L}(x)$ with $2\in \mathcal{L}(x)$ is of the form $\{2\}$.  Thus $\max \mathcal{L}(x)=2$ which yields $\rho_2(R)=2$ and (3) holds.
Clearly (3) implies (4).  Assume (4) holds.  If $R$ has class number greater than 2, then Lemma \ref{carlitzlemma} implies that $\rho_2(R)\geq 3$.
It easily follows from Lemma \ref{rhok} item (4) and induction that $\rho_k(R)>k$ for all $k\geq 2$, a contradiction.  Thus $R$ has class number at most 2 and (1) holds.  Hence (1), (2), (3), and (4) are equivalent.  

If (3) holds, then $2\in \mathcal{L}(x)$ implies that $2=\max\mathcal{L}(x)$ and $\mathcal{L}(x)=\{2\}$, which yields (5).   Statements (5) and (6) are equivalent by the definition of the length set.  If (6) holds, then $|\mathcal{L}(xy)|=1$ implies
that $\max \mathcal{L}(xy)=2$, which in turn yields (3).  This completes the proof.
\end{proof}

Let's take a slightly different look at the length set.  Given an algebraic number ring $R$ and $x$ a nonzero nonunit, suppose that
\[
\mathcal{L}(x)=\{n_1, \ldots ,n_k\}
\]
where $n_1<n_2 < \cdots < n_k$.  The delta set of $x$ is defined as 
\[
\Delta(x)=\{n_i - n_{i-1}\, |\, 2\leq i \leq k\}
\]
with $\Delta(x)=\emptyset$ if $k=1$.  We can convert this local descriptor into a global one by setting
\[
\Delta(R) =\bigcup_{x\in R^\bullet} \Delta(x).
\]
When $R$ is a Krull domain (a more general structure than an algebraic number ring) a great deal is known about the 
structure of $\Delta(R)$ (see \cite[Section 6.7]{GHKb} and \cite{Sch}).  

We show how the notion of the $\Delta$-set fits in with class number 2.    

\begin{theorem}
Let $R$ be an algebraic number ring.  Then $R$ has class number at most 2 if and only if $\Delta (R)=\emptyset$.
\end{theorem}

\begin{proof}   The implication ($\Rightarrow$) clearly holds by Carlitz's Theorem.  For ($\Leftarrow$), if $\Delta(R)=\emptyset$ and $R$ has class number greater than 2, 
then Lemma \ref{carlitzlemma} yields a contradiction.  This completes the proof.
\end{proof}

%%%%%%%%%%%%%%%%%%%% omega section %%%%%%%%%%%%%%%%%%%%%%%%%%%%%%%%%%%%%%%%%%%%

\section{Beyond the Length Set.}  Our characterizations to this point have been solely dependent
on the length set.  We now consider an invariant that relies on individual factorizations as much as or 
more than the set $\mathcal{L}(x)$.  It offers a numeric measure of how far an 
element is from being prime.  

\begin{definition} \label{firstdef} Let $R$ be an algebraic number ring.
For $x\in R^\bullet$, we define
$\omega(x) = n$ if $n$ is the smallest positive integer with the
property that whenever $x\mid a_1\cdots a_t$, where each
$a_i\in \mathcal{A}(R)$, there is a $T\subseteq \{1,2,\dots, t\}$ with
$|T|\leq n$ such that $x\mid \prod_{k\in T}a_k$. If no such $n$
exists, then $\omega(x)=\infty$.  For $x\in \mathcal{U}(R)$, we define
$\omega(x)=0$.  Finally, set
\[
\omega(R)=\sup\{\omega(\alpha)\, |\, \alpha\in \mathcal{A}(R)\}.
\]
\end{definition}

\noindent The definition above is taken from \cite{CSP}, but there are several other equivalent versions that can be found in the literature (see \cite{AndChap}).  It follows directly from the definition that an element $x\in R$ is prime if and only if $\omega(x)=1$.  The survey paper \cite{OnPel} is a good general reference on the $\omega$-function and we illustrate Definition \ref{firstdef} by appealing directly to the class number 2 case.

\begin{example}\label{omegaex}
Suppose that $R$ is an algebraic number ring of class number 2.  We use the classification of irreducible elements of $R$ given in Example \ref{ex1} to determine the $\omega$-values of the irreducibles of $R$.  If $\alpha$ is a prime element, then $\omega(\alpha)=1$.  So, let $\alpha$ be a nonprime element of $\mathcal{A}(R)$ where $(\alpha)=\mf{p}^2$ for a nonprincipal prime ideal $\mf{p}$ of $R$.  Thus $\omega(x)>1$, so suppose that $\alpha \, |\, \beta_1\cdots \beta_r$ where each $\beta_i$ is irreducible in $R$ and $r\geq 2$.  Hence, either one of the $\beta_i$'s is of the form $(\beta_i)=\mf{p}^2$, or there are irreducibles $\beta_i$ and $\beta_j$ (with $i\neq j$) so that $(\beta_i)=\mf{p}\mf{q}_1$ and $(\beta_j)=\mf{p}\mf{q}_2$ where $\mf{q}_1$ and $\mf{q}_2$ are nonprincipal prime ideals of $R$ distinct from $\mf{p}$.  In the first case, $\alpha$ is an associate of $\beta_i$ and in the second, $\alpha \, |\, \beta_i\beta_j$ and hence $\omega(\alpha)=2$.  A similar argument shows that $\omega(\alpha)=2$ if $(\alpha)=\mf{p}\mf{q}$ where $\mf{p}$ and $\mf{q}$ are distinct nonprincipal prime ideals of $R$.  
\end{example}

We introduce an aid
which will simplify the computation of $\omega(x)$.

\begin{definition}
Let $x\in R^\bullet$ where $R$ is an algebraic number ring.  
A \textit{bullet} for $x$ is a product $\beta_1\cdots \beta_r$ of irreducible elements $\beta_1, \ldots ,\beta_r$ of $R$ such that
\begin{enumerate}
\item[(i)] $x$ divides the product $\beta_1\cdots \beta_r$, and
\item[(ii)] for each $1\leq i\leq r$, $x$ does not divide $\beta_1\cdots \beta_r/\beta_i$.
\end{enumerate}
The set of bullets of $x$ is denoted $\mathrm{bul}(x)$.
\end{definition}

The notion of bullet gives us a nice tool to compute $\omega(x)$.   To see this, if $\beta_1\cdots \beta_r$ is a bullet for $x\in R^\bullet$, then $x$ divides no product of the form $\beta_1\cdots \beta_r/\beta_i$ for any $i$, and by definition
$\omega(x) \geq r$.  On the other hand, if $\alpha_1\cdots \alpha_t$ is a product of $t$ irreducibles of $R$ with 
$x\, |\, \alpha_1\cdots \alpha_t$ and $\alpha_1\cdots \alpha_t$ is not a bullet of $x$, then some subproduct of 
$\alpha_1\cdots \alpha_t$ must be a bullet.  We have essentially shown the following (a complete proof can be found in \cite[Proposition 2.10]{OnPel}) .

\begin{proposition}\label{hilo} If $R$ is an algebraic number ring and $x\in R^\bullet$, then 
\[
\omega(x) = \sup\{r\, |\, \beta_1\cdots \beta_r\in \mathrm{bul}(x)\mbox{  where each  }\beta_i\in \mathcal{A}(R)\}.
\]
\end{proposition}

Hence, for $R$ with class number 2, Example \ref{omegaex} shows that $\omega(R)=2$.  Proposition \ref{hilo} implies another nice
finiteness condition.

\begin{corollary}\label{almosthilo} Let $R$ be an algebraic number ring and $x\in\mathcal{A}(R)$.  Then
\[
\omega(x)\leq D(\mathcal{C}(R))<\infty
\]
and hence $\omega(R) \leq D(\mathcal{C}(R))<\infty$.  
\end{corollary}

\noindent In fact, the interested reader can find a proof that $\omega(R)=D(\mathcal{C}(R))$ in \cite[Corollary 3.3]{AndChap}.

\begin{proof}[Proof of Corollary \ref{almosthilo}] We prove only the first assertion, as the second follows directly from it.
Let $x\in \mathcal{A}(R)$.  
Write $(x)=\mf{p}_1^{t_1}\cdots \mf{p}_k^{t_k}$ for distinct prime ideals $\mf{p}_1, \ldots, \mf{p}_k$ in $R$. 
Let $\alpha_1, \ldots ,\alpha_n$ be irreducibles of $R$ such that $x\, |\, \alpha_1\cdots \alpha_n$.  For each $\mf{p}_i$ choose a minimal subset $T_i\subseteq \{1, \ldots, n\}$ so that $\mf{p}_i^{t_i}\, |\, (\prod_{j\in T_i}\alpha_j )$.  
Set $A_i=\{\alpha_j\, |\, j\in T_i\}$.  By the minimality of $T_i$, each $(\alpha_j)$, with $\alpha_j$ in $A_i$, is divisible by $\mf{p}_i$ and hence $|A_i|\leq t_i$.   If $A=\cup_{j=1}^k A_i$, then by Corollary \ref{thecor}, $|A|\leq t_1+\cdots +t_k \leq D(\mathcal{C}(R))$.
By using the multiplicative properties
of prime ideals, we obtain that $x\, |\, \prod_{\alpha_i\in A} \alpha_i$,       
which completes the proof.
\end{proof}

A slight adjustment in the proof of Corollary \ref{almosthilo} yields a class number 2 characterization (see \cite[Theorem 3.4]{AndChap}).

\begin{theorem}
Let $R$ be an algebraic number ring.  Then $R$ has class number at most 2 if and only if $\omega(R)\leq 2$.
\end{theorem}

\begin{proof}
While the argument is trivial using the remark directly following Corollary \ref{almosthilo}, for completeness we offer a proof.
Our work in Example \ref{omegaex}, along with the fact that class number 1 trivially implies $\omega(R)=1$, yields ($\Rightarrow$).  
For ($\Leftarrow$), assume $\omega(R)\leq 2$ and that $R$ has class number greater than 2.  We pivot in a manner similar to Lemma \ref{carlitzlemma}.  Suppose $\mathcal{C}(R)$ has an element $g$ with $|g|=n>2$.  Let 
$\mf{p}_1, \ldots ,\mf{p}_n$ be distinct prime ideals of $R$ with $[\mf{p}_i]=g$.  Let $x\in \mathcal{A}(R)$ be such
that $(x)=\mf{p}_1\cdots \mf{p}_n$.  If for each $1\leq i\leq n$ the irreducible $\alpha_i$ is such that 
$(\alpha_i)=\mf{p}_i^n$, then it is clear that $\alpha_1\cdots \alpha_n$ is a bullet for $x$ and $\omega(x)\geq n>2$.  
If $\mathcal{C}(R)$ only has nontrivial elements of order 2, then let $\alpha, \beta_1, \beta_2$, and $\beta_3$ be the irreducibles constructed in the second part of the proof of Lemma \ref{carlitzlemma}. 
As in the previous case, $\beta_1\beta_2\beta_3$ is a bullet for $\alpha$ and $\omega(x)\geq 3>2$.
In either case, $\omega(R) >2$, which completes the proof.
\end{proof}

\section{The Grand Finale!}  In these pages we have accomplished a lot.  To demonstrate, we tie it all together 
in one last tribute to class number 2.

\begin{Nutshell}
Let $R$ be an algebraic number ring.  The following statements are equivalent.
\begin{enumerate}
\item[(1)] $R$ has class number at most 2.
\item[(2)] $R$ is a half-factorial domain.
\item[(3)] $\rho(R)=1$.
\item[(4)] $\rho_2(R) = 2$.
\item[(5)] $\rho_k(R) = k$ for some $k\geq 2$.
\item[(6)] For all irreducibles $x$ and $y$ in $R$, $\mathcal{L}(xy) = \{2\}$.
\item[(7)] For all irreducibles $x$ and $y$ in $R$, $|\mathcal{L}(xy)|=1$.
\item[(8)] $\Delta (R)=\emptyset$.
\item[(9)]  $\omega(R)\leq 2$.
\end{enumerate}
\end{Nutshell}

We note that our work has not endeavored to determine exactly how many irreducible factorizations there are of an element $x$ in a class number 2 algebraic number ring $R$.  If $R$ has class number 2, then a formula 
for this computation is contained in \cite{CHR}.  
A detailed study of the asymptotic behavior of factorizations in rings with class number 2 can be found in \cite{HK}.
A more general approach to counting irreducible factorizations (with no restrictions on the class number) can be found in \cite{mystery}.

\begin{acknowledgment}{Acknowledgments.} The author gratefully acknowledges support under an Academic Leave during the fall of 2017 funded by Sam Houston State University.
He would also like to thank the referees and editor Susan Colley
for comments that greatly improved the exposition of this paper.
\end{acknowledgment}

\begin{biog}
\item[Scott Chapman] is Scholar in Residence and Distinguished Professor of Mathematics
at Sam Houston State University in Huntsville, Texas. In December of 2016 he finished
a five year appointment as Editor of the American Mathematical Monthly. His editorial
work, numerous publications in the area of non-unique factorizations, and years of
directing REU Programs, led to his designation in 2017 as a Fellow of the American
Mathematical Society.
\begin{affil}
Department of Mathematics and Statistics, Sam Houston State University, Box 2206, Huntsville, TX  77341\\
Scott.chapman@shsu.edu
\end{affil}
\end{biog}
\vfill\eject

\end{document}